\documentclass[a4paper,10pt]{amsart}
\usepackage[all]{xy}
\usepackage{amssymb}
\newtheorem{definition}{Definition}[section]
\newtheorem{lemma}[definition]{Lemma}
\newtheorem{prop}[definition]{Proposition}
\newtheorem{theorem}[definition]{Theorem}
\newtheorem{cor}[definition]{Corollary}

\newtheorem{prob}[definition]{Problem}
\newtheorem{fact}[definition]{Fact}

\newcommand{\class}{\mathop{{\rm class}}\nolimits}
\newcommand{\con}{\mathop{{\rm con}}\nolimits}
\newcommand{\id}{\mathop{{\rm id}}\nolimits}
\newcommand{\Jac}{\mathop{{\rm Jac}}\nolimits}
\renewcommand{\k}{\mathop{\raisebox{0pt}{\rm k}}\nolimits}

\newcommand{\Tr}{\mathop{{\rm Tr}}\nolimits}

\newcommand{\ov}{\overline}
\newcommand{\sep}{\supseteq}
\newcommand{\seq}{\subseteq}

\newcommand{\xr}{\xrightarrow}

\newcommand{\Ra}{\Rightarrow}

\newcommand{\al}{\alpha}
\newcommand{\be}{\beta}

\newcommand{\Del}{\Delta}
\newcommand{\eps}{\varepsilon}
\newcommand{\ga}{\gamma}

\newcommand{\om}{\omega}
\renewcommand{\phi}{\varphi}

\newcommand{\Sig}{\Sigma}

\newcommand{\C}{\mathbb{C}}
\newcommand{\Q}{\mathbb{Q}}
\newcommand{\Z}{\mathbb{Z}}
\DeclareMathOperator{\Mat}{M}

\newcommand{\lt}{\left}
\newcommand{\rt}{\right}

\begin{document}

\footskip30pt

\date{Thu Dec 15 12:30:36 EST 2005}

\title{Idempotent ideals and non-finitely generated projective
modules over integral group rings of polycyclic-by-finite groups}

\author{Peter~A. Linnell}

\address{Department of Mathematics, Virginia Tech, Blacksburg,
VA 24061-0123, USA} \email{linnell@math.vt.edu}

\author{Gena Puninski}

\address{School of Mathematics, The University of Manchester, Lamb Building,
Booth Road East, Manchester M13 9PL. England}
\email{gpuninski@maths.man.ac.uk}

\author{Patrick Smith}

\address{Department of Mathematics, University of Glasgow, UK}
\email{pfs@maths.gla.ac.uk}

\thanks{The second author is partially supported by the NSF grant no.~0501224.}


\begin{abstract}
We prove that every non-finitely generated projective module over
the integral group ring of a polycyclic-by-finite group $G$ is
free if and only if $G$ is polycyclic.
\end{abstract}

\maketitle

\pagestyle{plain}

\section{Introduction}\label{Int}

A question about existence of idempotent ideals is often
embroidered in a definition of certain classes of rings. For
example (see \cite[Prop.~6.3]{M-R}), a hereditary noetherian prime
ring $R$ is Dedekind prime, if $R$ has no nontrivial idempotent
ideals. However, it is usually difficult to decide whether a ring
in a given class of rings has a nontrivial idempotent ideal.

For instance, Akasaki \cite[Thm.~2.1]{Aka1} proved that the
integral group ring of any finite non-soluble group has a
nontrivial idempotent ideal. This result was completed by
Roggenkamp \cite{Rog}: the integral group ring of a finite soluble
group has no nontrivial idempotent ideals. Later P.~Smith
\cite[Thm.~2.2]{Smi} extended Roggenkamp's result to integral
group rings of polycyclic groups. In fact, using arguments similar
to Akasaki's, it is not difficult to conclude (see
Proposition~\ref{char}) that the integral group ring of a
polycyclic-by-finite group $G$ has no nontrivial idempotent ideals
if and only if $G$ is polycyclic (i.e., if $G$ is soluble).

By Bass \cite{Bass}, it is quite often that non-finitely generated
projective modules over a given noetherian ring are free. For
instance (see \cite[Cor.~4.5]{Bass}), this is true for modules
over an indecomposable commutative noetherian ring. Confirming
this observation, Swan \cite[Thm.~7]{Swan} proved that every
non-finitely generated projective module over the integral group
ring of a finite soluble group is free. Unfortunately, contrary to
what was expected (see Bass \cite[p.~24]{Bass}), some infinitely
generated non-free modules over integral group rings of
non-soluble finite groups were found by Akasaki \cite{Aka2} and
Linnell \cite{Lin}. Namely, they proved that, if $G$ is a finite
non-soluble group, then there exists an infinitely generated
projective $\Z G$-module which is not free.

Moreover, Akasaki \cite{Aka2} revealed the role played by
idempotent ideals in those kind of results. The main idea was to
use an idempotent ideal $I$ of $\Z G$ and a remarkable theorem of
Whitehead \cite{Whi} to guarantee the existence of a projective
$\Z G$-module $P$ whose trace is equal to $I$. The proof was
completed by the following observation: every finitely generated
projective module over the integral group ring of a finite group
is a generator (hence $P$ has no finitely generated direct
summands).

In this paper we extend this Akasaki and Linnell result to the
class of (infinite) polycyclic-by-finite groups. Namely (see
Theorem~\ref{last}), we prove that every non-finitely generated
projective module over the integral group ring of a
polycyclic-by-finite group $G$ is free if and only if $G$ is
polycyclic (that is, if $G$ is soluble). Moreover, if $G$ is not
polycyclic, then there exists a projective $\Z G$-module without
finitely generated direct summands.

One implication of this theorem (that every non-finitely generated
projective module over the integral group ring of a polycyclic
group is free) was proved by the third author \cite[Thm.~2.3]{Smi}.
However the proof was a little brief, so we have included a full
proof of this result using \cite[Lemma 4.6]{Smi1}, which is a
criterion for a projective module to be finitely generated.

We also give another proof of this implication based on Cohn's
universal localizations. Though it is longer, it may be used to
extend the results of the paper to the case of soluble groups.

To prove the other implication, extending Akasaki's result we show
(in Proposition~\ref{every}) that every finitely generated
projective module over the integral group ring of a
polycyclic-by-finite group is a generator. As a key preliminary
step, using $K$-theory, we prove (in Corollary~\ref{st-free})
that, if $G$ is a polycyclic-by-finite group and $P$ is a finitely
generated projective $\Z G$-module, then $(P\otimes_{\Z} \C)^n$ is
a free $\C G$-module for some $n$.

Despite the existence of `bad' projective modules over some
integral group rings (say, over $\Z A_5$) being established, the
structure of these modules is widely unknown. For instance, the
old question of P.~Linnell (Problem~8.34 in \cite{Kou}) about the
existence of a non-finitely generated indecomposable projective
module over the integral group ring of a finite group is still
open.

To make the paper self-contained and available even to neophytes,
we included most basic definitions and explanations from group
theory and the theory of projectives modules.

\section{Projective modules}

All rings will be associative rings with unity $1$, and most
modules will be (unitary) right modules.

Recall that a module $F$ over a ring $R$ is \emph{free}, if it is
isomorphic to a direct sum of copies of $R$, that is, $F\cong
R^{(\al)}$ for some cardinal $\al$. A module $P$ is said to be
\emph{projective}, if $P$ is isomorphic to a direct summand of a
free module. For instance, if $f\in R$ is an idempotent, then
$R_R=fR\oplus (1-f)R$, hence $fR$ is a projective (right)
$R$-module. Every finitely generated projective module $P$ over a
ring $R$ is isomorphic to the module $E R^n$, where $E$ is an
idempotent $n\times n$ matrix over $R$, and $R^n$ is a column
of height $n$ over $R$.

The trace of $P$, $\Tr(P)$, is defined as the sum of images of
all morphisms from $P$ to $R_R$. For instance, if $f\in R$ is an
idempotent and $P=fR$, then $\Tr(P)=RfR$ is the two-sided ideal of
$R$ generated by $f$. More generally, if a finitely generated
projective module $P$ is given by an idempotent matrix $E$ as
above, then $\Tr(P)$ is the two-sided ideal of $R$ generated by
the entries of $E$.

\begin{fact}(see \cite[Prop.~2.40]{Lam})\label{PI}
If $P$ is a projective $R$-module, then $I=\Tr(P)$ is an
idempotent ideal. Moreover, $P=PI$, and $I$ is the least two-sided
ideal of $R$ with this property.
\end{fact}

Thus, the trace of any projective module is an idempotent ideal.
However, the question whether a given idempotent ideal is the
trace of some projective module is hard to answer. For example, if
$R$ is a commutative ring and $I$ is a finitely generated
idempotent ideal of $R$, then (see \cite[Lemma~2.43]{Lam}) $I$ is
generated by an idempotent $f$, hence $I=\Tr(fR)$.

If $R$ is not commutative, and $I$ is an idempotent ideal of $R$
finitely generated on one side, it is not always true that $I$ is
generated by an idempotent. But such an $I$ is a trace of a
projective module as the following fact shows.

\begin{fact}\cite[Cor.~2.2]{Whi}\label{Whi}
Let $I$ be a two-sided idempotent ideal of a ring $R$ such that $I$ is
finitely generated as a left ideal. Then there exists a projective
right $R$-module $P$ such that $\Tr(P)=I$.
\end{fact}

Following Bass, we say that a projective module $P$ over a ring
$R$ is \emph{$\om$-big} (uniformly $\om$-big in the terminology of
Bass) if, for every two-sided ideal $I$ of $R$, the (projective)
$R/I$-module $P/PI$ is not finitely generated. In particular, $P$
itself is not finitely generated.

If $R$ is a ring, $\Jac(R)$ will denote the Jacobson radical of
$R$. The following fact shows that $\om$-big projective modules
are often free.

\begin{fact}\cite[Thm.~3.1]{Bass}\label{big}
Let $R$ be a ring such that $R/\Jac(R)$ is right noetherian. Then
every $\om$-big countably generated projective right $R$-module is
free.
\end{fact}

A nonzero (projective) module $P$ over a ring $R$ is said to be
\emph{stably free} (of rank $m$), if $P\oplus R^n\cong R^{n+m}$
for some $n, m$. If $R$ is right noetherian, then $m$ is the ratio
of Goldie dimensions of $P$ and $R_R$.

We say that $P$ is a \emph{generator} if, for some $k$, there is
an epimorphism $P^k\to R_R$. Clearly this is the same as
$\Tr(P)=R$. It is easily checked (as T.~Stafford pointed out to
the second author) that every stably free module over a (right)
noetherian ring is a generator. We need the following refinement
of this fact.

\begin{fact}\cite[Thm.~11.1.3]{M-R}\label{sfree}
Let $P$ be a stably free module over a right noetherian ring of
Krull dimension $h$. If the rank of $P$ exceeds $h$, then $P$ is
free.
\end{fact}

The following theorem of Kaplansky imposes severe restrictions on
the size of indecomposable projective modules.

\begin{fact}\cite[Cor.~2.6.2]{A-F}\label{sum}
Every projective module is a direct sum of countably generated
modules.
\end{fact}

A ring $R$ is said to be \emph{local}, if $R$ has a unique maximal
right (and left) ideal, that is, if $R/\Jac(R)$ is a skew field.

\begin{fact}\cite[Cor.~26.7]{A-F}\label{local}
Every projective module over a local ring is free.
\end{fact}

We also need the following result on projective modules.

\begin{fact}\cite[Lemma~4]{LLS}\label{fg}
Let $R$ be a subring of the ring $S$ and let $P$ be a projective
$R$-module. If the induced (projective) $S$-module $P\otimes_{R}
S$ is finitely generated, then $P$ is finitely generated.
\end{fact}

\section{Groups}

We will use $e$ for the unity of a group $G$, and then $\{e\}$
will denote the subgroup of $G$ consisting of $e$.

A series of groups $\{e\}=G_0\subset G_1\subset \dots \subset G_n$
is said to be a \emph{subnormal series} of a group $G$, if $G_n=G$,
and each $G_i$ is a normal subgroup of $G_{i+1}$, $i=0, \dots, n-1$.

A group $G$ is \emph{soluble}, if it has a subnormal series with
abelian factors $G_{i+1}/G_i$.  The \emph{commutant} $G'=[G,G]$ of
a group $G$ is defined as the (normal) subgroup of $G$ generated by
all commutators $[a,b]=aba^{-1}b^{-1}$, $a, b\in G$. Analogously,
setting $G^{(n)}= (G^{(n-1)})'$, one can define the \emph{$n$th
commutant} of $G$, $G^{(n)}$. It is easily checked that each
$G^{(n)}$ is a normal subgroup of $G$.

\begin{fact}\label{sol}
A group $G$ is soluble if and only if $G^{(n)}=\{e\}$ for some
$n$.
\end{fact}

A group $G$ is said to be \emph{perfect}, if $G=[G,G]$, that is,
$G$ coincides with its commutant. For instance, every alternating
group $A_n$, $n\geq 5$ is perfect. Thus, if $G$ is a finite
non-soluble group, then some commutant of $G$ is a nonidentity perfect
normal subgroup.

A group $G$ is called \emph{polycyclic}, if $G$ has a
subnormal series with cyclic factors. Since every cyclic group is
abelian, every polycyclic group is soluble. If all factors in a
subnormal series of $G$ are either cyclic or finite, $G$ is said
to be a \emph{polycyclic-by-finite} group. For instance, every
polycyclic group is polycyclic-by-finite. The number of infinite
cyclic factors in a subnormal series of a polycyclic-by-finite group
$G$ is called \emph{the Hirsch number} of $G$, $h(G)$.

\begin{fact}\label{poly-sol}
A polycyclic-by-finite group is soluble if and only if it is
polycyclic.
\end{fact}

If $G$ is a polycyclic-by-finite group, then the factors in a
subnormal series of $G$ can be rearranged such that all of them
but the last are cyclic. Thus every polycyclic-by-finite group $G$
has a normal polycyclic subgroup of finite index.

Now we extend a standard characterization of non-soluble finite groups
to general soluble-by-finite groups. This characterization should be
known, but we were not able to find a precise reference.

\begin{lemma}\label{nsol-char}
Let $G$ be a soluble-by-finite group which is not soluble.
Then $G$ has a nonidentity perfect normal subgroup.
\end{lemma}
\begin{proof}
First recall that, if $H$, $K$ and $L$ are normal subgroups of (any) group
$G$, then $[HK,L]=[H,L]\cdot [K,L]$. Indeed, the inclusion $\sep$ is clear,
and $\seq$ is easily checked using the following commutator identity:
$[hk,l]= h[k,l]h^{-1}\cdot [h,l]$.

Now the claim of the lemma is trivial, if $G$ is abelian, and also true
(see after Fact~\ref{sol}), if $G$ is finite.

Otherwise, let $G_1$ be a normal soluble subgroup of finite index in $G$,
and let $A$ be the last nonidentity term in the derived series of $G_1$. Thus
$A$ is a nontrivial normal abelian subgroup of $G$. Arguing by induction
(on the length of the derived series of $G_1$), we may assume that the result
has been already proved for $G/A$, that is, $G/A$ has a nonidentity perfect normal
subgroup $H/A$.

Thus $A\subset H$ are normal subgroups of $G$ such that $H$ is perfect modulo
$A$, hence $H=H^{(n)}A$ for every $n$. We prove that $H'$ is a perfect normal
subgroup of $G$. Then $H=H'A$ would imply that $H'\neq \{e\}$.

Indeed, since $H=H^{(2)}A$, by the above remark we obtain
$[H,A]=[H^{(2)}A,A]=[H^{(2)},A]\cdot [A,A]=[H^{(2)},A]\seq H^{(2)}$.
Then, using the same remark,
$H'=[H,H]=[H'A,H'A]=[H',H'A]\cdot [A,H'A]=H^{(2)} [H',A]\cdot [A,H']\cdot [A,A]=
H^{(2)}\cdot [H',A]\seq H^{(2)}$, hence $H'$ is perfect.
\end{proof}

It is known that every polycyclic-by-finite group $G$
is residually finite, but we need a stronger version of this
result. A group $G$ is said to be \emph{conjugacy separable}, if
for any $g, h\in G$ such that $g$ is not a conjugate of $h$ in
$G$, there is a normal subgroup $H$ of finite index in $G$ such
that $\bar g$ is not a conjugate of $\bar h$ in $G/H$.

\begin{fact}\cite{Rem}\label{Rem}
Every polycyclic-by-finite group $G$ is conjugacy separable.
Moreover, if $\ga_1, \dots, \ga_k$ are different conjugacy classes
of $G$, then there exists a normal subgroup $H$ of finite index in
$G$ such that $\bar \ga_1, \dots, \bar \ga_k$ are different
conjugacy classes in $G/H$, where $\bar{\gamma}_i$ indicates the
image of $\gamma_i$ in $G/H$.
\end{fact}

\section{Integral group rings}

If $G$ is a group, then $\Z G$ will denote the integral group ring
of $G$. Thus, every element $r\in \Z G$ has a unique
representation $r= \sum_{i=1}^k n_i g_i$, where $n_i\in \Z$ and
$g_i\in G$. In particular $1\cdot e$ is the unity of $\Z G$ and we
will usually write $1$ instead. There is a natural epimorphism
(the augmentation map) $\eps \colon \Z G\to \Z$ given by $\eps(r)=
\sum_{i=1}^k n_i$. The kernel of this map is the
\emph{augmentation ideal} $\Del(G)$ of $G$. Thus $r\in \Del(G)$ if
and only if $\sum_{i=1}^k n_i =0$. It is easily seen that
$\Del(G)$ is the free abelian group generated by
$\{g-1\mid g\in G\}$.  Furthermore $\eps$ extends to an epimorphism
of $n\times n$ matrices $\Mat_n(\mathbb{Z}G) \to \Mat_n(\mathbb{Z})$
by applying $\eps$ to each matrix entry.

If $H$ is a subgroup of $G$, the inclusion $H\seq G$ induces an
inclusion $\Z H\seq \Z G$. Let $\om(H)$ denote the two-sided ideal
of $\Z G$ generated by $\Del(H)$. If $H$ is normal in $G$, it is
easily checked that $\om(H)=\Z G\cdot \Del(H)= \Del(H)\cdot \Z G$.
A natural epimorphism $G\to G/H$ induces an epimorphism of rings
$\pi \colon \Z G \to \Z (G/H)$. The following fact describes the
kernel of $\pi$.

\begin{fact}(see \cite[Lemma~1.8]{Pas0})\label{kernel}
If $H$ is a normal subgroup of a group $G$, then we have the
following exact sequence: $0\to \om(H)\to \Z G \xr{\pi} \Z
(G/H)\to 0$. In particular, if $H=G$, we obtain $0\to \Del(G)\to
\Z G \xr{\eps} \Z\to 0$.
\end{fact}

We recall some basic facts on integral group rings of
polycyclic-by-finite groups.

\begin{fact}\cite[Thm.~5.12]{M-R}\label{noet}
If $G$ is a polycyclic-by-finite group, then $\Z G$ is a
noe\-the\-ri\-an ring.
\end{fact}

We also need the following powerful result of Roseblade.

\begin{fact}\cite[Cor.~C3]{Ros}\label{max}
Let $G$ be a polycyclic-by-finite group, and let $M$ be a maximal
two-sided ideal of $\Z G$. Then there exists a prime $p$ and a
normal subgroup $H$ of finite index in $G$ such that $p\,\Z$,
$\om(H)\seq M$.
\end{fact}

The proof of the following lemma is similar to Akasaki
\cite[Thm.~2.1]{Aka1}.

\begin{lemma}\label{idemp}
Let $G$ be a polycyclic-by-finite group which is not polycyclic.
Then the integral group ring of $G$ has a nontrivial idempotent
ideal.
\end{lemma}
\begin{proof}
Since $G$ is not polycyclic, it is not soluble (see
Fact~\ref{poly-sol}). Thus, by Fact~\ref{nsol-char}, $G$ has a nonidentity normal perfect
subgroup $H$, that is, $H=[H,H]$.

First we check that the augmentation ideal $\Del(H)$ of $\Z H$ is
idempotent. Indeed, as we have already noticed, the elements
$h-1$, $h\in H$ generate $\Del(H)$ as an abelian group. Thus, it
suffices to check that $h-1\in \Del(H)^2$. If
$h=[a,b]=aba^{-1}b^{-1}$ for some $a, b\in H$, then
$h-1=ab(a^{-1}b^{-1}-b^{-1}a^{-1})$ Since
$cd-dc=(c-1)(d-1)-(d-1)(c-1)\in \Del(H)^2$ for all $c, d\in H$, we
conclude that $h-1\in \Del(H)^2$.

Suppose that $h$ is an arbitrary element of $H$. Then
$h=h_1\cdot\ldots\cdot h_k$, where $h_i=[a_i,b_i]$ for some $a_i,
b_i\in H$. Using an induction on $k$, we may assume that $k\geq 2$
and $h'-1\in \Del(H)^2$, where $h'=h_2\cdot\ldots\cdot h_k\in H$.
Then $h-1=(h_1-1)(h'-1)+(h_1-1)+(h'-1)\in \Del(H)^2$ by what we
have already proved.

Now, as we have previously noticed, $\om(H)=\Z G\cdot \Del(H)=
\Del(H)\cdot \Z G$. Since $\Del(H)^2=\Del(H)$, one easily derives
that $\om(H)$ is a nonzero idempotent two-sided ideal of
$\Z G$.  From the obvious inclusion $\om(H)\seq \Del(G)$ it
follows that $\om(H)\neq \Z G$.
\end{proof}

Now we are in a position to characterize integral group rings of
polycyclic-by-finite groups possessing a nontrivial idempotent
ideal.

\begin{prop}\label{char}
Let $G$ be a polycyclic-by-finite group. Then the following are
equivalent:

1) $\Z G$ has no nontrivial idempotent two-sided ideals.

2) $G$ is polycyclic.
\end{prop}
\begin{proof}
1) $\Ra$ 2) follows from Lemma~\ref{idemp}, and 2) $\Ra$ 1) is
\cite[Thm.~2.2.]{Smi}.
\end{proof}

\section{One implication}

Recall that $K_0(R)$ denotes the Grothendieck group of finitely
generated projective modules over a ring $R$, that is, an abelian
group whose generators are finitely generated projective
$R$-modules subject to relations $[P]=[Q]+[Q']$ for all short
exact (hence split) sequences $0\to Q\to P\to Q'\to 0$. The
Grothendieck group $G_0(R)$ of $R$ is a similar abelian group
defined on the set of finitely generated $R$-modules.

Note that (see \cite[Cor.~5.60]{Lam}), if $R$ is a two-sided
noetherian ring, then the left and right global dimensions of $R$
coincide and can be calculated as the supremum of projective
dimensions of cyclic $R$-modules. If this supremum is finite, $R$
is said to be a \emph{regular} ring.

\begin{fact}\cite[Prop.~3.8]{Pas}\label{reg}
Let $R$ be a noetherian regular ring. Then the natural map $c
\colon K_0(R)\to G_0(R)$ sending $[P]$ to $[P]$ is an isomorphism
called the \emph{Cartan map}.
\end{fact}

In general the Cartan map will not be mono or epi. What is
important for us is the following instance of Fact~\ref{reg}.

\begin{fact}\cite[Thm.~3.13]{Pas0}\label{poly-reg}
If $G$ is a polycyclic-by-finite group, then the ring $\C G$ is
regular. Thus $K_0(\C G)$ is canonically isomorphic to $G_0(\C
G)$.
\end{fact}

If $R\seq S$ are rings and $P$ is a finitely generated $R$-module,
it induces a finitely generated $S$-module $P\otimes_R S$. If $S$
is flat as a left $R$-module, this gives rise to the induction map
$G_0(R)\to G_0(S)$, and similarly for $K_0$. For instance, if $H$
is a subgroup of a group $G$, and $S$ is a ring, then induction
via the inclusion $SH\seq SG$ induces maps $G_0(SH)\to G_0(SG)$
and $K_0(SH)\to K_0(SG)$.

The following fact is known as J.~Moody's induction theorem.

\begin{fact}(see \cite[Thm.~3.6.9]{Pas})\label{Moody}
Let $S$ be a noetherian ring and let $G$ be a polycyclic-by-finite
group. Suppose that $G_1, \dots, G_k$ are representatives of the
conjugacy classes of the maximal finite subgroups of $G$. Then
$G_0(SG)$ is generated by the images of various $G_0(SG_i)$ under
the induction map.
\end{fact}

If $R$ is a ring, then $H_0(R)=R/[R,R]$ will denote the factor of
$R$ by the additive subgroup generated by additive commutators
$[r,s]=rs-sr$ for $r,s \in R$. If $R=SG$ and $S$ is commutative,
then $[R,R]$ is an $S$-submodule of $R$, and $H_0(R)$ is
isomorphic to the free $S$-module $\class_0(G)$ generated by the set
$\con(G)$ of conjugacy classes of $G$.

The Hattori--Stallings map, $HS$, sends finitely generated
projective $R$-modules to elements of $H_0(R)$ in the following
way. Let $P$ be a finitely generated projective $R$-module, hence
$P\cong E R^n$ for some idempotent $n\times n$ matrix $E$ over
$R$. Then the trace of $E$ is an element $r\in R=SG$, and $HS(P)$
is the image $\bar r$ of $r$ in $H_0(R)=R/[R,R]$ (see \cite{Sta}
for more on that). For example, if $S$ is a field, then $HS(P)$
is an element of the $S$-vector space $\class_0(G)$ with $\con(G)$
as a basis.

The Hattori--Stallings map induces a homomorphism of abelian groups
$HS \colon K_0(R)\to H_0(R)$ (note that, even if $S$ is commutative,
$K_0(R)$ carries no natural $S$-module structure). For instance,
the image of $[R]\in K_0(R)$ is $\bar 1 = 1\cdot \{e\}$.

\begin{lemma}\label{Luck}
If $G$ is a polycyclic-by-finite group, then the map $K_0(\C
G)\otimes_{\Z} \C\to \class_0(G)$ induced by the
Hattori--Stallings map is an embedding.
\end{lemma}
\begin{proof}
By \cite[Lemma~2.15]{Luck} we have the following commutative diagram:

$$
\vcenter{%
\xymatrix@C=16pt@R=24pt{%
\varinjlim K_0(\C H)\otimes_{\Z} \C
\ar[r]^(.55){f}\ar[d]^{\cong}_{h}&
K_0(\C G)\otimes_{\Z} \C\ar[d]^{HS}\\
\class_0(G)_f\ar[r]^{g}&\class_0(G)\\
}}
$$

\vspace{2mm}

Here the direct limit is taken over the set of all finite
subgroups $H$ of $G$ with respect to inclusion and $G$-conjugation
(see \cite[p.~91]{B-L}). The morphism $f$ is induced by induction
maps $K_0(\C H)\to K_0(\C G)$. Furthermore, $\class_0(G)_f$
denotes the $\C$ vector space spanned by the conjugacy classes of
$G$ whose elements have finite order, and $g$ is a natural
inclusion.

As is explained in \cite[Thm.~2.22]{Luck} (or alternatively use
\cite[Prop.~of \S 1.8]{B-L}), J.~Moody's induction
theorem (Fact \ref{Moody}) and Fact~\ref{poly-reg} imply that $f$
is an isomorphism. Since $g$ is an embedding, the same is true for
$HS$.
\end{proof}

By Swan \cite[Thm~8.1]{Swan0}, if $P$ is a finitely generated
projective module over the integral group ring of a finite group
$G$, then $P\otimes_{\Z } \Q$ is a free $\Q\, G$-module. Thus the
following is a weak form of this result.

\begin{lemma}\label{st-free}
Let $P$ be a finitely generated projective module over the
integral group ring of a polycyclic-by-finite group $G$. Then
$(P\otimes_{\Z} \C)^n$ is a stably free $\C G$-module for some
$n$.
\end{lemma}
\begin{proof}
Let $P'=P\otimes_{\Z} \C$. Clearly $P'$ is a finitely generated
projective $\C G$-module, that is $[P']\in K_0(\C G)$. The kernel
of the natural map $K_0(\C G)\to K_0(\C G)\otimes_{\Z} \C$ is the
torsion part of $K_0(\C G)$. Thus, by Lemma~\ref{Luck}, it
suffices to prove that $HS(P') \in \mathbb{Z}\bar 1$.

Let $HS(P')=\sum_{i=1}^k \al_i \ga_i$, where $0\neq \al_i\in \C$
and $\ga_i$ are different conjugacy classes of $G$. If $HS(P')
\notin \mathbb{Z}\bar 1$,
then there is a conjugacy class $\gamma_l \ne \{e\}$
(where $1\le l \le k$). By Fact~\ref{Rem}, there is a normal
subgroup $H$ of finite index in $G$ such that $\bar{\ga}_1, \dots,
\bar{\ga}_n$ are different conjugacy classes of $G/H$, in
particular $\bar \ga_l\neq \{\bar e\}$.

Now $P'/P'\om(H)$ is a finitely generated projective $\C
(G/H)$-module (see Fact~\ref{kernel}), and clearly
$P'/P'\om(H)\cong P/P\om(H)\otimes_{\Z} \C$. Furthermore,
$P/P\om(H)$ is a finitely generated projective $\Z (G/H)$-module.
Then, by the aforementioned Swan's result, $P/P\om(H)\otimes_{\Z}
\Q$ is a free $\Q\, [G/H]$-module, hence $P'/P'\om(H)$ is a free
$\C(G/H)$-module. But
$HS(P'/P'\om(H))=\sum_{i=1}^k \al_i \bar \ga_i$ has a nonzero
coefficient by $\bar \ga_l\neq \{\bar e\}$, a contradiction.
\end{proof}

\begin{cor}\label{free}
Let $P$ be a finitely generated projective module over the
integral group ring of a polycyclic-by-finite group $G$. Then
$(P\otimes_{\Z} \C)^n$ is a free $\C G$-module for some positive
integer $n$.
\end{cor}
\begin{proof}
By Lemma~\ref{st-free}, $Q=(P\otimes_{\Z} \C)^n$ is a stably free
$\C G$-module for some $n$, and we may assume that $n\geq mh+1$,
where $h$ is the Hirsch number of $G$ and $m$ is the Goldie
dimension of $\C G$.

Then the Goldie dimension of $Q$ is at least $mh+1$, hence the
rank of $Q$ (which is the ratio of Goldie dimensions of $Q$ and $\C
G$) exceeds $h$. Also, by \cite[Prop.~6.1.1]{M-R}, the Krull
dimension of $\C G$ is equal to $h$. It remains to apply
Fact~\ref{sfree}.
\end{proof}

By Akasaki \cite[Cor.~1.4]{Aka2} every finitely generated
projective module over the integral group ring of a finite group
is a generator. We generalize this result as follows.

\begin{prop}\label{every}
Let $P$ be a finitely generated projective module over the
integral group ring of a polycyclic-by-finite group $G$. Then $P$
is a generator.
\end{prop}
\begin{proof}
Otherwise $I=\Tr(P)$ is a nontrivial two-sided ideal of $\Z G$
such that $P=PI$. If $M$ is any maximal ideal of $\Z G$ containing
$\Tr(P)$, then $P=PM$. By Fact~\ref{max}, there is a normal
subgroup $H$ of finite index in $G$ such that $\om(H)\seq M$.
Write $\overline{M} = M/\omega(H)$, the image of $M$ in
$\mathbb{Z}G/\omega(H) \cong \mathbb{Z}(G/H)$. Then $P'=P/P\om(H)$
is a finitely generated projective $\Z (G/H)$-module such that
$P'\ov{M}=P'$. But $G/H$ is a finite group, hence (by Akasaki's
result) $P'=0$ and we deduce that $(P\otimes_{\mathbb{Z}}
\mathbb{C})^n \omega(H) = (P \otimes_{\mathbb{Z}} \mathbb {C})^n$
for all positive integers $n$.  Using Corollary~\ref{free}, we
conclude that $P \otimes_{\mathbb{Z}} \mathbb{C} = 0$,
consequently $P = 0$ which is a contradiction.
\end{proof}

Now we prove the first implication of the foregoing theorem.

\begin{theorem}\label{main}
Let $G$ be a polycyclic-by-finite group which is not polycyclic.
Then there exists a projective $\Z G$-module $P$ such that $P$ has
no finitely generated direct summands. In particular, $P$ is not
free.
\end{theorem}
\begin{proof}
By (the proof of) Lemma~\ref{idemp}, there exists a nontrivial
normal subgroup $H$ of $G$ such that $I=\om(H)$ is an idempotent
ideal. Since $\Z G$ is noetherian, by Fact~\ref{Whi} there exists
a projective $\Z G$-module $P$ such that $\Tr(P)=I$.

Suppose that $Q$ is a finitely generated direct summand of $P$. By
Proposition~\ref{every}, $\Tr(Q)=\Z G$, hence $\Tr(P)=\Z G$, a
contradiction.
\end{proof}

\section{The other implication}

In this section we prove the following proposition.

\begin{prop}\label{one-imp}
Let $P$ be a non-finitely generated projective module over the
integral group ring of a polycyclic group $G$. Then $P$ is free.
\end{prop}

But first we recall some definitions and results of \cite{Smi1}.
Let $I$ be a two-sided ideal of a ring $R$. We define a descending
chain of two-sided ideals $I^{\al}$ as follows: put
$I^{\al+1}=I^{\al}\cdot I$, and $I^{\al}=\bigcap_{\be<\al} I^{\be}$,
if $\al$ is a limit ordinal. If $\ga$ is the least ordinal such that
$I^{\ga}=I^{\ga +1}$, then denote $I^{\ga}$ by $\k^1(I)$.

Note that this definition is right handed: we obtain a left
variant of it by setting $I^{\al+1}=I\cdot I^{\al}$ in the
above definition. For every positive integer $n$ define
$\k^{n+1}(I)=\k^1(\k^n(I))$.

Let $\id(I)$ denote the sum of idempotent ideals of $R$ contained
in $I$. Then $\id(I)$ is the unique maximal idempotent ideal
contained in $I$ and $\id(I)\seq \k^n(I)$ for all $n$. A ring
$R$ is called \emph{right shallow}, if for each ideal $I$ of
$R$ there exists a positive integer $n$ such that
$\k^n(I)=\id(I)$.

\begin{fact}\cite[Cor.~2.10]{Smi1}\label{shallow}
If $G$ is a finite group, then the integral group ring of $G$
is right and left shallow.
\end{fact}

\begin{cor}\label{I-sol}
Let $G$ be a finite soluble group and let $I$ be a proper two-sided
ideal of $\Z G$. Then there exists a positive integer $n$ such
that $\k^n(I)=0$.
\end{cor}
\begin{proof}
By Fact~\ref{shallow}, $\k^n(I)=\id(I)$ for some $n$. But, by
Roggenkamp~\cite[Prop.]{Rog} the integral group ring of a finite
soluble group has no nontrivial idempotent ideals.
\end{proof}

Before proving the next lemma we need the following instance of
B.~Hartley's result.

\begin{fact}\cite[Thm.~E]{Har}\label{Del-ab}
If $G$ is a finitely generated abelian group and $\Del=\Del(\Z G)$,
then $\bigcap_n \Del^n=0$, in particular $\k^1(\Del)=0$.
\end{fact}

\begin{lemma}\label{I-poly}
Let $G$ be a polycyclic group. Then for every proper ideal $I$
of $\Z G$ there exists a positive integer $n$ such that
$\k^n(I)=0$.
\end{lemma}
\begin{proof}
Let $M$ be a maximal ideal of $\Z G$ such that $I\seq M$. By
Fact~\ref{max}, there exists a normal subgroup $H$ of finite
index in $G$ such that $\om(H)\seq M$. Then
$\bar I= I+\om(H)/\om(H)$ is a proper ideal of
$\Z G/\om(H)\cong \Z(G/H)$. By Corollary~\ref{I-sol}, there exists
a positive integer $m$ such that $\k^m(\bar I)=\bar 0$, that is,
$\k^m(I)\seq \om(H)$. Thus it suffices to prove that
$\k^l(\om(H))=0$ for some $l$. Moreover, since
$\om(H)= \Z G\cdot \Del(H)= \Del(H)\cdot \Z G$, it is enough to
show that $\k^l(\Del(H))=0$.

We prove this by induction on the derived length of $H$. If $H$
is abelian, the result follows from Fact~\ref{Del-ab}. Otherwise
$H'\neq \{e\}$ and we may assume that $\k^s(\Del(H'))=0$ for some
$s$, hence $\k^s(\om(H'))=0$. Since $H/H'$ is a finitely generated
abelian group, $\k^1(\Del(H/H'))=0$ by Fact~\ref{Del-ab}. Thus
$\k^1(\Del(H))\seq \om(H')$, therefore
$\k^{s+1}(\Del(H))\seq \k^s(\om(H'))=0$.
\end{proof}

Now we are ready to prove Proposition~\ref{one-imp}.

By Fact~\ref{sum} we may assume that $P$ is countably generated.
Moreover, by Fact~\ref{big} it suffices to prove that $P$ is
$\om$-big, that is, $P/PI$ is not finitely generated for every
proper two-sided ideal $I$ of $\Z G$.

Indeed, if $P/PI$ is finitely generated then, by
\cite[Lemma~4.6]{Smi1}, $P/P\k^n(I)$ is finitely generated for every
$n$. But, by Lemma~\ref{I-poly}, $\k^n(I)=0$ for some $n$, hence
$P$ is finitely generated, a contradiction.

Now we are in a position to formulate the main result of the
paper.

\begin{theorem}\label{last}
Let $G$ be a polycyclic-by-finite group. Then the following are
equivalent:

1) $G$ is not polycyclic;

2) there exists a non-finitely generated projective $\Z G$-module
that is not free;

3) there exists a projective $\Z G$-module $P$ such that $P$ has
no finitely generated direct summands.
\end{theorem}
\begin{proof}
By Theorem~\ref{main} and Proposition~\ref{one-imp}.
\end{proof}

\section{Localizations}

Let $I$ be a prime ideal of a noetherian ring $R$. Then, by
Goldie's theorem, $R/I$ has a classical quotient ring $Q=Q(R/I)$
which is a simple artinian ring. Let $\Sig$ be the set of all
square $R$-matrices that are regular modulo $I$ (that is, become
invertible over $Q$). Then Cohn's universal localization
$R_{\Sig}$ can be defined as a ring homomorphism $R\to R_{\Sig}$
such that 1) the images of all matrices in $\Sig$ are invertible
over $R_{\Sig}$, and 2) for every ring homomorphism $R\to S$ which
inverts all matrices in $\Sig$, there exists a unique ring
homomorphism $R_{\Sig}\to S$ completing the following commutative
diagram:

$$
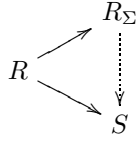
\vcenter{%
\xymatrix@C=20pt@R=6pt{%
&R_{\Sig}\ar@{.>}[dd]&\\
R\ar[ur]\ar[dr]\\
& S\\
}}
$$

\vspace{2mm}

The ring $R_{\Sig}$ is unique (as a universal object) and can be
constructed by adding formal inverses to matrices in $\Sig$. For a
more constructive way of building $R_{\Sig}$ see Malcolmson
\cite{Mal}.

\begin{fact}\label{Cohn}\cite[Thm.~4.1]{Cohn}
The map $R\to R_{\Sig}$ induces an embedding $R/I\to
R_{\Sig}/\Jac(R_{\Sig})$ such that $R_{\Sig}/\Jac(R_{\Sig})\cong
Q(R/I)$.
\end{fact}

Note that, if the prime ideal $I$ is (right) localizable, then
$R_{\Sig}\cong R_I$, the localization of $R$ at $I$.

If $M$ is an $R$-module, we define the localization of $M$,
$M_{\Sig}$, setting $M_{\Sig}=M\otimes_R R_{\Sig}$. Thus
$M_{\Sig}$ is an $R_{\Sig}$-module.

A group $G$ is said to be \emph{poly-$\Z$}, if it has a subnormal
series with factors isomorphic to $\Z$. Thus, every poly-$\Z$
group is polycyclic. Furthermore, every polycyclic group contains
a normal poly-$\Z$ subgroup of finite index.

\begin{fact}(see \cite[Lemma~3.7.8]{Pas})\label{Ore}
Every poly-$\Z$ group $G$ is torsion-free and the integral group
ring of $G$ is an Ore domain.
\end{fact}

\begin{lemma}\label{mono}
Let $H$ be a poly-$\Z$ group, and let $f_A \colon \Z H^n\to \Z
H^n$ be a homomorphism of (free) $\Z H$-modules given by left
multiplication by an $n\times n$ matrix $A$ over $\Z H$. If $f_A$
is a monomorphism modulo $\Del(H)$, then $f_A$ is a monomorphism.
\end{lemma}
\begin{proof}
By Strebel \cite[Prop.~1.5]{Str}, $H$ is in the class $D(\Z)$,
that is, any map $f \colon P\to Q$ between projective $\Z
H$-modules is a monomorphism, if the induced map of abelian groups
$f\otimes 1 \colon P\otimes_{\Z H} \Z \to Q\otimes_{\Z H} \Z$ is a
monomorphism.

It remains to notice that $f_A\otimes 1$ is an endomorphism of a
$\Z$-module $\Z^n\cong \lt[\Z H/\Del(H)\rt]^n$ given by left
multiplication by $\eps(A)$.
\end{proof}

We put in use the following result.

\begin{fact}\cite[Thm.~7]{Swan}\label{Swan}
Every non-finitely generated projective module over the integral
group ring of a finite soluble group  is free.
\end{fact}

Now we give a different proof of Proposition~\ref{one-imp}.

As above we may assume that $P$ is countably generated,
and we have to show that $P/PI$ is not finitely generated for
every proper two-sided ideal $I$ of $\Z G$.

Suppose that $P/PI$ is finitely generated. If $M$ is a maximal
two-sided ideal of $\Z G$ containing $I$, then $P/PM$ is also
finitely generated. We prove that this leads to a contradiction.

By Fact~\ref{max}, there is a prime $p$ and a normal subgroup $H$
of finite index in $G$ such that $p\, \Z$, $\om(H)\seq M$.
Moreover (see the remark above), we may assume that $H$ is a
poly-$\Z$ group.

Note that $P/P \om(H)$ is a projective $\Z
G/\om(H)=\Z(G/H)$-module and $G/H$ is a finite soluble group. If
$P/P \om(H)$ is not finitely generated, it is a free infinite rank
$\Z (G/H)$-module by Fact~\ref{Swan}. Since $\om(H)\seq M$, we
conclude that $P/PM$ is a free infinite rank $\Z G/M$-module, a
contradiction.

Thus $P/P\om(H)$ is a finitely generated $\Z G$-module. Since $H$
is of finite index in $G$ and $P\omega(H) = P\Delta(H)$, we deduce
that $P/P\Del(H)$ is a finitely generated $\Z H$-module.

Set $J=p\,\Z + \Del(H)$. Then $\Z H/ J\cong \Z/p\,\Z$, hence $J$
is a (completely) prime ideal of $\Z H$. Since $P/P\Del(H)$ is
finitely generated, $P/PJ$ is a finite dimensional vector space
over $\Z/p\,\Z$ and hence $P/PJ$ is finite.

Let $\Sig$ denote the set of all square $\Z H$-matrices that are
invertible modulo $J$. By Fact~\ref{Cohn}, the universal
localization $\Z H_{\Sig}$ is a local ring such that $\Z H_{\Sig}/
\Jac(\Z H_{\Sig})\cong \Z/p\, \Z$.

Now $P_{\Sig}=P\otimes_{\Z H} \Z H_{\Sig}$ is a projective module
over the local ring $\Z H_{\Sig}$, hence $P_{\Sig}$ is free by
Fact~\ref{local}.  Also the image of $P$ in $P_{\Sig}$ generates
$P_{\Sig}$ as a $\Z H_{\Sig}$-module.  Since $P/PJ$ is finite and
the image of $PJ$ in $P_{\Sig}/P_{\Sig} \Jac(P_{\Sig})$ is 0
(because the image of $J$ in $\Z H_{\Sig}/\Jac(Z H_{\Sig})$ is 0;
cf.~Fact~\ref{Cohn}), we deduce that
$P_{\Sig}/P_{\Sig}\Jac(P_{\Sig})$ is a finitely generated $\Z
H_{\Sig}$-module. Therefore $P_{\Sig}$ is a finitely generated
free $\Z H_{\Sig}$-module.

Below we prove that the natural morphism $\Z H\to \Z H_{\Sig}$ is
an embedding. Then Fact~\ref{fg} implies that $P$ is a finitely
generated $\Z H$-module, the desired contradiction.

By Fact~\ref{Ore}, $\Z H$ is an Ore domain. Let $F$ be the
classical ring of quotients of $\Z H$ which is a skew field. Since
$\Z H$ is an Ore domain, we can identify it with the universal
localization with respect to the set $\Sig'$ of all square $\Z
H$-matrices which are left (=right) regular, that is, whose images
in $F$ are invertible.

We show that $\Sig\seq \Sig'$. Suppose that $A$ is an $n\times n$
matrix in $\Sig$. Then left multiplication by $A$ induces a
morphism of free $\Z H$-modules $f_A \colon \Z H^n \to \Z H^n$
that is a monomorphism modulo $J$. Clearly this map is a
monomorphism modulo $\Del(H)$ (since every square $\Z$-matrix is
regular, if it is regular modulo $p\, \Z$).  Then, by
Lemma~\ref{mono}, $f_A$ is a monomorphism, hence $A\in \Sig'$.

Thus $\Sig\seq \Sig'$, hence, by the universal property, there
exists a morphism of rings $\Z H_{\Sig}\to \Z H_{\Sig'}=F$
completing the following commutative diagram:

$$
\vcenter{%
\xymatrix@C=14pt@R=6pt{%
&\Z H_{\Sig}\ar@{.>}[dd]&\\
\Z H\ar[ur]\ar[dr]\\
& \Z H_{\Sig'}\\
}}
$$

\vspace{2mm}

Since $\Z H$ is a subring of $F$, the morphism $\Z H\to \Z
H_{\Sig}$ is an embedding, as desired.

We wonder whether the results of this paper can be extended
from polycyclic-by-finite to arbitrary soluble groups.

\begin{prob}\label{conj}
Is every non-finitely generated projective module over the
integral group ring of a soluble group free?
\end{prob}

\end{document}